\documentclass[10pt,twoside,reqno]{amsart}
\usepackage{amssymb}
\usepackage{multirow}
\calclayout

\numberwithin{equation}{section}
\makeatletter

\renewcommand{\@secnumfont}{\bfseries}

\renewcommand{\section}{\@startsection{section}{1}%
  {0mm}{.7\linespacing\@plus\linespacing}{.5\linespacing}
  {\normalfont\bfseries\centering}}

\newcommand{\bibsection}{\@startsection{section}{1}%
  {0mm}{.7\linespacing\@plus\linespacing}{.5\linespacing}
  {\normalfont\scshape\centering}}

\renewcommand{\@biblabel}[1]{#1.}

\newtheorem{thm}{\bf Theorem}[section]

\usepackage{cite}

\begin{document}

\vspace{1.3cm}

\title{Degenerate Changhee numbers and polynomials of the second kind}

\author{Taekyun Kim}
\address{Department of Mathematics, Kwangwoon University, Seoul 139-701, Republic
	of Korea}
\email{tkkim@kw.ac.kr}

\author{Dae San Kim}
\address{Department of Mathematics, Sogang University, Seoul 121-742, Republic
	of Korea}
\email{dskim@sogang.ac.kr}

\begin{abstract}
In this paper, we consider the degenerate Changhee numbers and polynomials of the second kind which are different from the previously introduced degenerate Changhee numbers and polynomials by Kwon-Kim-Seo (see [11]). We investigate some interesting identities and properties for these numbers and polynomials. In addition, we give some new relations between the degenerate Changhee polynomials of the second kind and the Carlitz's degenerate Euler polynomials.
\end{abstract}

\subjclass[2010]{11B83; 11S80}
\keywords{degenerate Changhee polynomials of the second kind,  }
\maketitle

\markboth{\centerline{\scriptsize Degenerate Changhee numbers and polynomials of the second kind}}
{\centerline{\scriptsize T. Kim, D. S. Kim}}

\section{Introduction}
Let $p$ be a fixed odd prime number. Throughout this paper, $\mathbb{Z}_p$, $\mathbb{Q}_p$ and $\mathbb{C}_p$ will denote the ring of $p$-adic integers, the field of $p$-adic rational numbers and the completion of an algebraic closure of $\mathbb{Q}_p$. The $p$-adic norm $|\cdot |_p$ is normalized by $|p|_p = \frac{1}{p}$. Let $C(\mathbb{Z}_p)$ be the space of continuous functions on $\mathbb{Z}_p$. For $f \in C(\mathbb{Z}_p)$, the fermionic $p$-adic integral on $\mathbb{Z}_p$ is defined by Kim as 

\begin{equation}\begin{split}\label{01}
I(f) &= \int_{\mathbb{Z}_p} f(x)   d\mu_{-1} (x) = \lim_{N \rightarrow \infty} \sum_{x=0}^{p^N-1} f(x) \mu_{-1} (x+p^N \mathbb{Z}_p)\\
&=\lim_{N \rightarrow \infty} \sum_{x=0}^{p^N-1} f(x) (-1)^x,\quad (\textnormal{see}\;\; [8,19]                                                                                                                                                                                                                                                                             ).
\end{split}\end{equation}

From \eqref{01}, we note that

\begin{equation}\begin{split}\label{02}
I(f_n)+(-1)^{n-1}I(f) = 2 \sum_{a=0}^{n-1} (-1)^{n-1-a} f(a),\quad (\textnormal{see} \,\, [8,18,19]),
\end{split}\end{equation}
where $f_n(x) = f(x+n)$, $(n \in \mathbb{N})$. It is well known that the Euler polynomials are defined by the generating function 
\begin{equation}\begin{split}\label{03}
\frac{2}{e^t+1} e^{xt} = \sum_{n=0}^\infty E_n(x) \frac{t^n}{n!},\quad (\textnormal{see} \,\, [ 1-20]).
\end{split}\end{equation}
When $x=0$, $E_n=E_n(0)$ are called the Euler numbers. 

In [2,3], L. Carlitz considered the degenerate Euler polynomials given by the generating function 

\begin{equation}\begin{split}\label{04}
\frac{2}{(1+\lambda t)^{\frac{1}{\lambda }}+1}(1+\lambda t)^{\frac{x}{\lambda }} = \sum_{n=0}^\infty \mathcal{E}_{n,\lambda }(x) \frac{t^n}{n!},\,\,(\lambda  \in \mathbb{R}).
\end{split}\end{equation}

When $x=0$, $\mathcal{E}_{n,\lambda } = \mathcal{E}_{n,\lambda }(0)$ are called the degenerate Euler numbers. From \eqref{04}, we easily note that

\begin{equation*}\begin{split}
\sum_{n=0}^\infty \lim_{\lambda  \rightarrow 0} \mathcal{E}_{n,\lambda }(x) \frac{t^n}{n!}&=
\lim_{\lambda  \rightarrow 0} \frac{2}{(1+\lambda t)^{\frac{1}{\lambda }}+1}(1+\lambda t)^{\frac{x}{\lambda }} \\
&= \frac{2}{e^t+1}e^{xt} = \sum_{n=0}^\infty E_n(x) \frac{t^n}{n!}.
\end{split}\end{equation*}

Thus, we have

\begin{equation*}\begin{split}
\lim_{\lambda  \rightarrow 0} \mathcal{E}_{n,\lambda }(x) = E_n(x),\,\,(n \geq 0),\quad (\textnormal{see} \,\, [2]).
\end{split}\end{equation*}

As is well known, the Changhee polynomials are defined by the generating function 
\begin{equation}\begin{split}\label{05}
\frac{2}{t+2}(1+t)^x = \sum_{n=0}^\infty Ch_n(x) \frac{t^n}{n!},\quad (\textnormal{see} \,\, [ 7,9]).
\end{split}\end{equation}

When $x=0$, $Ch_n =Ch_n(0)$, $(n \geq 0)$, are called the Changhee numbers. From \eqref{02}, we note that

\begin{equation}\begin{split}\label{06}
\int_{\mathbb{Z}_p} (1+t)^{x+y}   d\mu_{-1} (y) = \frac{2}{2+t}(1+t)^x = \sum_{n=0}^\infty Ch_n(x) \frac{t^n}{n!}.
\end{split}\end{equation}

Thus, by \eqref{06}, we get
\begin{equation}\begin{split}\label{07}
\int_{\mathbb{Z}_p} (x+y)_n   d\mu_{-1} (y) = Ch_n(x), \,\,(n \geq 0), \quad (\textnormal{see} \,\, [7]),
\end{split}\end{equation}
where $(x)_0=1$, $(x)_n = x(x-1)\cdots(x-n+1)$, $(n \geq 1)$.

It is not difficult to show that
\begin{equation}\begin{split}\label{08}
\int_{\mathbb{Z}_p} e^{(x+y)t}   d\mu_{-1} (y) = \frac{2}{e^t+1} e^{xt} = \sum_{n=0}^\infty E_n(x) \frac{t^n}{n!}.
\end{split}\end{equation}
By \eqref{08}, we get
\begin{equation}\begin{split}\label{09}
\int_{\mathbb{Z}_p}  (x+y)^n  d\mu_{-1} (y)= E_n(x),\,\,(n \geq 0),\quad (\textnormal{see} \,\, [6,7,8,18,19]).
\end{split}\end{equation}

The Stirling numbers of the first kind are defined by
\begin{equation}\begin{split}\label{10}
(x)_n = \sum_{l=0}^n S_1(n,l) x^l, \,\,(n \geq 0),\quad (\textnormal{see} \,\, [1-20]),
\end{split}\end{equation}

and those of the second kind are given by

\begin{equation}\begin{split}\label{11}
x^n = \sum_{l=0}^n S_2(n,l) (x)_l,\,\,(n \geq 0),\quad (\textnormal{see} \,\, [1,4,5,7,18]).
\end{split}\end{equation}

From \eqref{06} and \eqref{08}, we note that

\begin{equation}\begin{split}\label{12}
Ch_n(x) = \sum_{l=0}^n E_l(x) S_1(n,l),
\end{split}\end{equation}
and
\begin{equation*}\begin{split}
E_n(x) = \sum_{l=0}^n Ch_n(x) S_2(n,l),\,\,(n \geq 0),\quad (\textnormal{see} \,\, [7]).
\end{split}\end{equation*}

Recently, the degenerate Changhee polynomials are introduced by Kwon-Kim-Seo as

\begin{equation}\begin{split}\label{13}
\frac{2\lambda }{2\lambda +\log(1+\lambda t)} \big(1+\log(1+\lambda t)^{\frac{1}{\lambda }}\big)^x = \sum_{n=0}^\infty Ch_{n,\lambda }^* (x) \frac{t^n}{n!}.
\end{split}\end{equation}

When $x=0$, $Ch_{n,\lambda }^* = Ch_{n,\lambda }^*(0)$ are called the degenerate Changhee numbers (see [11]).

Note that $\lim_{\lambda  \rightarrow 0} Ch_{n,\lambda }^* (x) = Ch_n(x)$, $(n \geq 0)$.

Recently, many researchers have studied Changhee numbers and polynomials (see [1-20]). In this paper, we consider the degenerate Changhee numbers and polynomials of the second kind which are different from the previous introduced degenerate Changhee numbers and polynomials by Kwon-Kim-Seo (see [11]). We give some new and interesting identities and properties for these numbers and polynomials. In addition, we investigate some relations between the degenerate Changhee polynomials of the second kind and Carlitz's degenerate Euler polynomials.

\section{Degenerate Changhee numbers and polynomials of the second kind}

From \eqref{02}, we note that

\begin{equation}\begin{split}\label{14}
\int_{\mathbb{Z}_p} \big(1+\lambda \log(1+t)\big)^{\frac{x+y}{\lambda }}  d\mu_{-1} (y) = \frac{2}{1+\big(1+\lambda \log(1+t)\big)^{\frac{1}{\lambda }}} \big(1+\lambda \log(1+t)\big)^{\frac{x}{\lambda }},
\end{split}\end{equation}
where $\lambda  \in \mathbb{C}_p$ with $|\lambda |_p \leq 1$. Now, we define the degenerate Changhee polynomials of the second kind by

\begin{equation}\begin{split}\label{15}
\frac{2}{1+\big(1+\lambda \log(1+t)\big)^{\frac{1}{\lambda }}} \big(1+\lambda \log(1+t)\big)^{\frac{x}{\lambda }} = \sum_{n=0}^\infty Ch_{n,\lambda }(x) \frac{t^n}{n!}.
\end{split}\end{equation}
Thus, by \eqref{14}, we get

\begin{equation}\begin{split}\label{16}
\int_{\mathbb{Z}_p} \big(1+\lambda \log(1+t)\big)^{\frac{x+y}{\lambda }}  d\mu_{-1} (y)  &=
\sum_{l=0}^\infty \int_{\mathbb{Z}_p}  {\frac{x+y}{\lambda } \choose l}  d\mu_{-1} (y) \lambda^l \big(\log(1+t)\big)^l \\
&=\sum_{n=0}^\infty \left( \sum_{l=0}^n S_1(n,l) \int_{\mathbb{Z}_p} {\frac{x+y}{\lambda } \choose l}l!d\mu_{-1} (y) \lambda ^l \right) \frac{t^n}{n!}.
\end{split}\end{equation}

Comparing the coefficients on both sides of \eqref{15} and \eqref{16}, we have

\begin{equation}\begin{split}\label{17}
\sum_{l=0}^n S_1(n,l) \int_{\mathbb{Z}_p} {\frac{x+y}{\lambda } \choose l}l! d\mu_{-1} (y) \lambda ^l  = Ch_{n,\lambda }(x),\,\,(n \geq 0).
\end{split}\end{equation}

The $\lambda $-analogue of falling factorial sequence is given by 
\begin{equation}\begin{split}\label{18}
(x)_{n,\lambda } = x(x-\lambda )\cdots(x-(n-1)\lambda ),\,\,(n \geq 1),\,\,(x)_{0,\lambda }=1.
\end{split}\end{equation}

Thus, by \eqref{17} and \eqref{18}, we get

\begin{equation}\begin{split}\label{19}
\sum_{l=0}^n S_1(n,l) \int_{\mathbb{Z}_p}    (x+y)_{l,\lambda } d\mu_{-1} (y) = Ch_{n,\lambda }(x),\,\,(n \geq 0).
\end{split}\end{equation}

\begin{thm}
For $n \geq 0$, we have
\begin{equation*}\begin{split}
\sum_{l=0}^n \int_{\mathbb{Z}_p}    (x+y)_{l,\lambda } d\mu_{-1} (y)  S_1(n,l)= Ch_{n,\lambda }(x).
\end{split}\end{equation*}
\end{thm}

It is not difficult to show that

\begin{equation}\begin{split}\label{20}
\int_{\mathbb{Z}_p} (1+\lambda t)^{\frac{x+y}{\lambda }}  d\mu_{-1} (y) &= \frac{2}{(1+\lambda t)^{\frac{1}{\lambda }}+1} (1+\lambda t)^{\frac{x}{\lambda }}\\
&= \sum_{n=0}^\infty \mathcal{E}_{n,\lambda }(x) \frac{t^n}{n!}.
\end{split}\end{equation}

From \eqref{20}, we note that
\begin{equation}\begin{split}\label{21}
\int_{\mathbb{Z}_p}  (x+y)_{n,\lambda }  d\mu_{-1} (y) = \mathcal{E}_{n,\lambda }(x),\,\,(n \geq 0).
\end{split}\end{equation}

Thus, from Theorem 2.1 and \eqref{21}, we obtain the following theorem.

\begin{thm}
For $n \geq 0$, we have
\begin{equation*}\begin{split}
\sum_{l=0}^n S_1(n,l) \mathcal{E}_{l,\lambda }(x) = Ch_{n,\lambda }(x).
\end{split}\end{equation*}
\end{thm}

By replacing $t$ by $e^t-1$ in \eqref{15}, we get
\begin{equation}\begin{split}\label{22}
\sum_{m=0}^\infty Ch_{m,\lambda }(x) \frac{1}{m!} (e^t-1)^m &= \frac{2}{(1+\lambda t)^{\frac{1}{\lambda }}+1}(1+\lambda t)^{\frac{x}{\lambda }} \\
&= \sum_{n=0}^\infty \mathcal{E}_{n,\lambda }(x) \frac{t^n}{n!}.
\end{split}\end{equation}

On the other hand,
\begin{equation}\begin{split}\label{23}
\sum_{m=0}^\infty Ch_{m,\lambda }(x) \frac{1}{m!} (e^t-1)^m &= \sum_{m=0}^\infty Ch_{m,\lambda }(x) \sum_{n=m}^\infty S_2(n,m) \frac{t^n}{n!}\\
&= \sum_{n=0}^\infty \left( \sum_{m=0}^n Ch_{m,\lambda }(x) S_2(n,m) \right) \frac{t^n}{n!}.
\end{split}\end{equation}

Therefore, by \eqref{22} and \eqref{23}, we obtain the following theorem.

\begin{thm}
For $n \geq 0$, we have
\begin{equation*}\begin{split}
\mathcal{E}_{n,\lambda }(x) = \sum_{m=0}^n Ch_{m,\lambda }(x) S_2(n,m).
\end{split}\end{equation*}
\end{thm}

When $x=0$, $Ch_{n,\lambda }(x) = Ch_{n,\lambda }(0)$, $(n \geq 0)$, are called the degenerate Changhee numbers of the second kind.

From \eqref{15}, we note that

\begin{equation}\begin{split}\label{24}
& \sum_{n=0}^\infty Ch_{n,\lambda }(x) \frac{t^n}{n!} = \frac{2}{1+\big(1+\lambda \log(1+t)\big)^{\frac{1}{\lambda }}} \big(1+\lambda \log(1+t)\big)^{\frac{x}{\lambda }}\\
&= \left( \sum_{l=0}^\infty Ch_{l,\lambda }\frac{t^l}{l!} \right)
\left( \sum_{m=0}^\infty {\frac{x}{\lambda }\choose m} \lambda ^m \big(\log(1+t)\big)^m \right)\\
&= \left( \sum_{l=0}^\infty Ch_{l,\lambda }\frac{t^l}{l!} \right)
\left( \sum_{m=0}^\infty (x)_{m,\lambda } \sum_{k=m}^\infty S_1(k,m) \frac{t^k}{k!} \right)\\
&= \left( \sum_{l=0}^\infty Ch_{l,\lambda }\frac{t^l}{l!} \right)
\left( \sum_{k=0}^\infty \left( \sum_{m=0}^k (x)_{m,k} S_1(k,m) \right) \frac{t^k}{k!} \right)\\
&= \sum_{n=0}^\infty \left( \sum_{k=0}^n \sum_{m=0}^k {n \choose k} (x)_{m,\lambda } S_1(k,m) Ch_{n-k,\lambda } \right) \frac{t^n}{n!}.
\end{split}\end{equation}

By comparing the coefficients on both sides of \eqref{24}, we obtain the following theorem.

\begin{thm}
For $n \geq 0$, we have
\begin{equation*}\begin{split}
Ch_{n,\lambda }(x) = \sum_{k=0}^n \sum_{m=0}^k {n \choose k} (x)_{m,\lambda } S_1(k,m) Ch_{n-k,\lambda }.
\end{split}\end{equation*}
\end{thm}

By \eqref{02} , we easily get
\begin{equation}\begin{split}\label{25}
\int_{\mathbb{Z}_p} f(x+1)   d\mu_{-1} (x) + \int_{\mathbb{Z}_p} f(x)   d\mu_{-1} (x) = 2f(0).
\end{split}\end{equation}

Thus, by \eqref{25}, we get

\begin{equation}\begin{split}\label{26}
\int_{\mathbb{Z}_p} \big(1+\lambda \log(1+t)\big)^{\frac{x+1}{\lambda }}  d\mu_{-1} (x) + \int_{\mathbb{Z}_p} \big(1+\lambda \log(1+t)\big)^{\frac{x}{\lambda }}  d\mu_{-1} (x) = 2. 
\end{split}\end{equation}

From \eqref{15} and \eqref{26}, we have
\begin{equation}\begin{split}\label{27}
\frac{2}{1+\big(1+\lambda \log(1+t)\big)^{\frac{1}{\lambda }}} \big(1+\lambda \log(1+t)\big)^{\frac{1}{\lambda }} + \frac{2}{1+\big(1+\lambda \log(1+t)\big)^{\frac{1}{\lambda }}}  = 2.
\end{split}\end{equation}

From \eqref{15} and \eqref{27}, we have
\begin{equation}\begin{split}\label{28}
\sum_{n=0}^\infty \big( Ch_{n,\lambda }(1)+Ch_{n,\lambda } \big) \frac{t^n}{n!} = 2.
\end{split}\end{equation}

Comparing the coefficients on both sides of \eqref{28}, we obtain the following theorem.

\begin{thm}
For $n \geq 0$, we have
\begin{equation*}\begin{split}
Ch_{n,\lambda }(1) + Ch_{n,\lambda } = \begin{cases}
2,&\text{if}\,\, n =0\\
0,&\text{if}\,\, n \geq 1.
\end{cases}
\end{split}\end{equation*}
\end{thm}

By Theorem 2.5, we easily get
\begin{equation*}\begin{split}
Ch_{0,\lambda }=1, \,\, Ch_{1,\lambda } = -\frac{1}{2},\,\, Ch_{2,\lambda }= \frac{1}{2} (1+\lambda ), \cdots.
\end{split}\end{equation*}

For $d \in \mathbb{N}$ with $d \equiv 1 $ (mod 2), by \eqref{02} we have
\begin{equation}\begin{split}\label{29}
\int_{\mathbb{Z}_p}  f(x+d)  d\mu_{-1} (x) + \int_{\mathbb{Z}_p}  f(x)  d\mu_{-1} (x) = 2 \sum_{a=0}^{d-1} (-1)^a f(a).
\end{split}\end{equation}

Let us take $f(x) =  \big(1+\lambda \log(1+t)\big)^{\frac{x}{\lambda }}$. Then by \eqref{29}, we get
\begin{equation}\begin{split}\label{30}
&\int_{\mathbb{Z}_p}  \big(1+\lambda \log(1+t)\big)^{\frac{x}{\lambda }}   d\mu_{-1} (x) \\
&=\frac{2}{ \big(1+\lambda \log(1+t)\big)^{\frac{d}{\lambda }}+1} \sum_{a=0}^{d-1} (-1)^a  \big(1+\lambda \log(1+t)\big)^{\frac{a}{\lambda }}\\
&= \sum_{a=0}^{d-1} (-1)^a \frac{2}{1+ \big(1+ \frac{\lambda }{d}d\log(1+t)\big)^{\frac{d}{\lambda }}} \big(1+ \tfrac{\lambda }{d}d\log(1+t)\big)^{\frac{d}{\lambda } \frac{a}{d}}.
\end{split}\end{equation}

By \eqref{20}, we easily get

\begin{equation}\begin{split}\label{31}
\frac{2}{1+ \big(1+ \frac{\lambda }{d}d\log(1+t)\big)^{\frac{d}{\lambda }}} \big(1+  \tfrac{\lambda }{d}d\log(1+t)\big)^{\frac{d}{\lambda } \frac{a}{d}}
&= \sum_{m=0}^\infty \mathcal{E}_{m, \frac{\lambda }{d}} \left( \frac{a}{d} \right) \frac{d^m}{m!} \big(\log(1+t)\big)^m\\
&= \sum_{m=0}^\infty \mathcal{E}_{m, \frac{\lambda }{d}} d^m \sum_{n=m}^\infty S_1(n,m) \frac{t^n}{n!}\\
&= \sum_{n=0}^\infty \left( \sum_{m=0}^n d^m  \mathcal{E}_{m, \frac{\lambda }{d}} S_1(n,m) \right) \frac{t^n}{n!}.
\end{split}\end{equation}

From \eqref{30} and \eqref{31}, we note that

\begin{equation}\begin{split}\label{32}
\sum_{n=0}^\infty Ch_{n,\lambda } \frac{t^n}{n!} &= \int_{\mathbb{Z}_p}  \big(1+\lambda \log(1+t)\big)^{\frac{x}{\lambda }}   d\mu_{-1} (x) \\
&= \sum_{n=0}^\infty \left( \sum_{m=0}^n d^m \sum_{a=0}^{d-1} (-1)^a \mathcal{E}_{m, \frac{\lambda }{d}} S_1(n,m) \right) \frac{t^n}{n!}.
\end{split}\end{equation}

Therefore, by \eqref{32}, we obtain the following theorem.

\begin{thm}
For $n \geq 0$, $d \in \mathbb{N}$ with $d \equiv 1$ (mod 2), we have
\begin{equation*}\begin{split}
Ch_{n,\lambda } = \sum_{m=0}^n d^m \sum_{a=0}^{d-1} (-1)^a \mathcal{E}_{m, \frac{\lambda }{d}} S_1(n,m).
\end{split}\end{equation*}
\end{thm}

Now, we observe that
\begin{equation}\begin{split}\label{33}
&\frac{2}{1+\big(1+\lambda \log(1+t)\big)^{\frac{1}{\lambda }}} \big(1+\lambda \log(1+t)\big)^{\frac{x+1}{\lambda }} + \frac{2\big(1+\lambda \log(1+t)\big)^{\frac{x}{\lambda }}}{1+\big(1+\lambda \log(1+t)\big)^{\frac{1}{\lambda }}} \\
&= 2\big(1+\lambda \log(1+t)\big)^{\frac{x}{\lambda }}.
\end{split}\end{equation}

Thus, by \eqref{15} and \eqref{33}, we get

\begin{equation}\begin{split}\label{34}
\sum_{n=0}^\infty \left( Ch_{n,\lambda }(x+1) + Ch_{n,\lambda }(x) \right) \frac{t^n}{n!} &=
2 \sum_{m=0}^\infty (x)_{m,\lambda } \frac{1}{m!} \big( \log(1+t)\big)^m\\
&= \sum_{n=0}^\infty \left( 2 \sum_{m=0}^n (x)_{m,\lambda } S_1(n,m) \right)  \frac{t^n}{n!}.
\end{split}\end{equation}

Thus, by \eqref{34}, we obtain the following theorem.

\begin{thm}
For $n \geq 0$, we have
\begin{equation*}\begin{split}
 Ch_{n,\lambda }(x+1) + Ch_{n,\lambda }(x) =  2 \sum_{m=0}^n (x)_{m,\lambda } S_1(n,m) .
\end{split}\end{equation*}
\end{thm}

From \eqref{29}, we have

\begin{equation}\begin{split}\label{35}
&\frac{2}{1+\big(1+\lambda \log(1+t)\big)^{\frac{1}{\lambda }}} \big(1+\lambda \log(1+t)\big)^{\frac{d}{\lambda }} + \frac{2}{1+\big(1+\lambda \log(1+t)\big)^{\frac{1}{\lambda }}}  
\\&= 2 \sum_{a=0}^{d-1} (-1)^a \big(1+\lambda \log(1+t)\big)^{\frac{a}{\lambda }},
\end{split}\end{equation}
where $d \in \mathbb{N}$ with $d \equiv 1$ (mod 2).

By \eqref{15} and \eqref{35}, we get

\begin{equation}\begin{split}\label{36}
\sum_{n=0}^\infty \left( Ch_{n,\lambda }(d) + Ch_{n,\lambda } \right) \frac{t^n}{n!} &= 
2 \sum_{a=0}^{d-1} (-1)^a \sum_{m=0}^\infty (a)_{m,\lambda } \frac{1}{m!} \big( \log(1+t) \big)^m\\ &= 2 \sum_{a=0}^{d-1} (-1)^a \sum_{n=0}^\infty \left( \sum_{m=0}^n (a)_{m,\lambda } S_1(n,m) \right) \frac{t^n}{n!}\\
&= \sum_{n=0}^\infty \left(  2 \sum_{a=0}^{d-1} \sum_{m=0}^n (a)_{m,\lambda } S_1(n,m)(-1)^a \right) \frac{t^n}{n!}.
\end{split}\end{equation}

Therefore, by \eqref{36}, we obtain the following theorem.

\begin{thm}
For $n \geq 0$, $d \in \mathbb{N}$ with $d \equiv 1$ (mod 2), we have
\begin{equation*}\begin{split}
Ch_{n,\lambda }(d) + Ch_{n,\lambda } = 2 \sum_{a=0}^{d-1} \sum_{m=0}^n (a)_{m,\lambda } S_1(n,m)(-1)^a.
\end{split}\end{equation*}
\end{thm}

Now, we consider the higher-order degenerate Changhee polynomials of the second kind which are derived from the multivariate fermionic $p$-adic integral on $\mathbb{Z}_p$.

For $r \in \mathbb{N}$, we define the higher-order degenerate Changhee polynomials of the second kind which are given by the multivariate fermionic $p$-adic integral on $\mathbb{Z}_p$ as follows:

\begin{equation}\begin{split}\label{37}
&\int_{\mathbb{Z}_p}\cdots \int_{\mathbb{Z}_p}  \big(1+\lambda \log(1+t)\big)^{\frac{x+x_1+\cdots+x_r}{\lambda }}   d\mu_{-1} (x_1)\cdots    d\mu_{-1} (x_r)\\
&= \left(\frac{2}{1+\big(1+\lambda \log(1+t)\big)^{\frac{1}{\lambda }}}\right)^r \big(1+\lambda \log(1+t)\big)^{\frac{x}{\lambda }} = \sum_{n=0}^\infty Ch_{n,\lambda }^{(r)}(x) \frac{t^n}{n!}.
\end{split}\end{equation}

When $x=0$, $Ch_{n,\lambda }^{(r)} = Ch_{n,\lambda }^{(r)}(0)$ are called the higher-order degenerate Changhee numbers of the second kind.

From \eqref{37}, we note that

\begin{equation}\begin{split}\label{38}
&\int_{\mathbb{Z}_p}\cdots \int_{\mathbb{Z}_p}  \big(1+\lambda \log(1+t)\big)^{\frac{x+x_1+\cdots+x_r}{\lambda }}   d\mu_{-1} (x_1)\cdots    d\mu_{-1} (x_r)\\
&= \sum_{m=0}^\infty \int_{\mathbb{Z}_p}\cdots \int_{\mathbb{Z}_p} {\frac{x_1+\cdots+x_r+x}{\lambda } \choose m}  d\mu_{-1} (x_1)\cdots    d\mu_{-1} (x_r) \lambda ^m \big(\log(1+t)\big)^m\\
&= \sum_{m=0}^\infty \int_{\mathbb{Z}_p}\cdots \int_{\mathbb{Z}_p} (x_1+\cdots+x_r+x)_{m,\lambda } d\mu_{-1} (x_1)\cdots    d\mu_{-1} (x_r) \frac{1}{m!} \big(\log(1+t)\big)^m\\
&= \sum_{n=0}^\infty \left( \sum_{m=0}^n  \int_{\mathbb{Z}_p}\cdots \int_{\mathbb{Z}_p} (x_1+\cdots+x_r+x)_{m,\lambda } d\mu_{-1} (x_1)\cdots    d\mu_{-1} (x_r) S_1(n,m)              \right) \frac{t^n}{n!}
\end{split}\end{equation}
It is easy to show that

\begin{equation}\begin{split}\label{39}
&\int_{\mathbb{Z}_p}\cdots \int_{\mathbb{Z}_p} (1+\lambda t)^{\frac{x_1+\cdots+x_r+x}{\lambda }} d\mu_{-1} (x_1)\cdots    d\mu_{-1} (x_r)\\
&= \left( \frac{2}{(1+\lambda t)^{\frac{1}{\lambda }}+1} \right)^r (1+\lambda t)^{\frac{x}{\lambda }} = \sum_{n=0}^\infty \mathcal{E}_{n,\lambda }^{(r)}(x) \frac{t^n}{n!},
\end{split}\end{equation}
where $\mathcal{E}_{n,\lambda }^{(r)}(x)$ are the Carlitz's degenerate Euler polynomials of order $r$.

Thus, by \eqref{39}, we get
\begin{equation}\begin{split}\label{40}
 \int_{\mathbb{Z}_p}\cdots \int_{\mathbb{Z}_p} (x_1+\cdots+x_r+x)_{m,\lambda } d\mu_{-1} (x_1)\cdots    d\mu_{-1} (x_r) = \mathcal{E}_{m,\lambda }^{(r)}(x),\,\,(m \geq 0).
\end{split}\end{equation}

Therefore, by \eqref{37}, \eqref{38} and \eqref{40}, we obtain the following theorem.

\begin{thm}
For $n \geq 0$, we have
\begin{equation*}\begin{split}
Ch_{n,\lambda }^{(r)}(x) = \sum_{m=0}^n \mathcal{E}_{m,\lambda }^{(r)}(x) S_1(n,m),\,\,(r \in \mathbb{N}).
\end{split}\end{equation*}
\end{thm}

By replacing $t$ by $e^t-1$ in \eqref{37}, we get
\begin{equation}\begin{split}\label{41}
&\int_{\mathbb{Z}_p}\cdots \int_{\mathbb{Z}_p} (1+\lambda t)^{\frac{x_1+\cdots+x_r+x}{\lambda }} d\mu_{-1} (x_1)\cdots    d\mu_{-1} (x_r) = \sum_{m=0}^\infty Ch_{m,\lambda }^{(r)}(x) \frac{(e^t-1)^m}{m!}\\
&= \sum_{n=0}^\infty \left( \sum_{m=0}^n Ch_{m,\lambda }^{(r)}(x) S_2(n,m) \right) \frac{t^n}{n!}.
\end{split}\end{equation}

Thus, by comparing the coefficients on both sides of \eqref{41} and \eqref{39}, we obtain the following theorem.

\begin{thm}
For $n \geq 0$, we have
\begin{equation*}\begin{split}
\mathcal{E}_{n,\lambda }^{(r)}(x) =  \sum_{m=0}^n Ch_{m,\lambda }^{(r)}(x) S_2(n,m).
\end{split}\end{equation*}
\end{thm}

The degenerate Stirling numbers of the second kind are defined by Kim as
\begin{equation}\begin{split}\label{43}
\frac{1}{n!} \big( (1+\lambda t)^{\frac{1}{\lambda }}-1 \big)^n = \sum_{m=n}^\infty S_{2,\lambda }(m,n) \frac{t^m}{m!},\quad (\textnormal{see} \,\, [10]).
\end{split}\end{equation}

Here the left hand side of \eqref{43} is given by

\begin{equation}\begin{split}\label{44}
&\frac{1}{n!} \big( (1+\lambda t)^{\frac{1}{\lambda }}-1 \big)^n  = \frac{1}{n!} \left( e^{\frac{1}{\lambda }\log(1+\lambda t)}-1 \right)^n\\
&= \sum_{l=n}^\infty S_2(l,n) \lambda ^{-l} \frac{1}{l!} \big(\log(1+\lambda t)\big)^l\\
&= \sum_{l=n}^\infty S_2(l,n) \lambda ^{-l} \sum_{m=l}^\infty S_1(m,l) \frac{\lambda ^m t^m}{m!}\\
&= \sum_{m=n}^\infty \left( \sum_{l=n}^m  S_2(l,n) \lambda ^{m-l} S_1(m,l) \right) \frac{t^m}{m!}.
\end{split}\end{equation}

Comparing \eqref{43} and \eqref{44}, we have

\begin{equation}\begin{split}\label{45}
S_{2,\lambda }(m,n) = \sum_{l=n}^m  S_2(l,n) \lambda ^{m-l} S_1(m,l),
\end{split}\end{equation}
where $m,n \geq 0$ with $m \geq n$.

Now, we observe that
\begin{equation}\begin{split}\label{46}
(1+\lambda t)^{\frac{x_1+\cdots+x_r+x}{\lambda }} &= \left( (1+\lambda t)^{\frac{1}{\lambda }}-1+1\right)^{x_1+\cdots+x_r+x}\\
&= \sum_{m=0}^\infty {x_1+\cdots+x_r+x \choose m} \big( (1+\lambda t)^{\frac{1}{\lambda }}-1\big)^m\\
&= \sum_{m=0}^\infty \left(x_1+\cdots+x_r+x \right)_m \sum_{n=m}^\infty S_{2,\lambda }(n,m) \frac{t^n}{n!}\\
&= \sum_{n=0}^\infty \left( \sum_{m=0}^n S_{2,\lambda }(n,m) \left(x_1+\cdots+x_r+x \right)_m \right) \frac{t^n}{n!}\\
&= \sum_{n=0}^\infty \left( \sum_{m=0}^n S_{2,\lambda }(n,m) (x_1+\cdots+x_r+x)_{m} \right) \frac{t^n}{n!}
\end{split}\end{equation}

Thus, by \eqref{41} and \eqref{46}, we get

\begin{equation}\begin{split}\label{47}
&\int_{\mathbb{Z}_p}\cdots \int_{\mathbb{Z}_p} (1+\lambda t)^{\frac{x_1+\cdots+x_r+x}{\lambda }} d\mu_{-1} (x_1)\cdots    d\mu_{-1} (x_r) \\
&= \sum_{n=0}^\infty \left( \sum_{m=0}^n S_{2,\lambda }(n,m) \int_{\mathbb{Z}_p}\cdots \int_{\mathbb{Z}_p} (x_1+\cdots+x_r+x)_{m} d\mu_{-1} (x_1)\cdots    d\mu_{-1} (x_r) \right) \frac{t^n}{n!}\\
&= \sum_{n=0}^\infty \left( \sum_{m=0}^n S_{2,\lambda }(n,m) Ch_{m}^{(r)} (x) \right) \frac{t^n}{n!}.
\end{split}\end{equation}

Therefore, by \eqref{41} and \eqref{47}, we obtain the following theorem.

\begin{thm}
For $n \geq 0$, we have
\begin{equation*}\begin{split}
\sum_{m=0}^n Ch_{m}^{(r)} (x) S_{2,\lambda }(n,m)= \sum_{m=0}^n Ch_{m,\lambda }^{(r)}(x) S_2(n,m).
\end{split}\end{equation*}
\end{thm}

From the generating function of the higher-order degenerate Changhee numbers of the second kind, we note that

\begin{equation}\begin{split}\label{48}
&\left( \frac{2}{\big(1+\lambda \log(1+t)\big)^{\frac{1}{\lambda }}+1}\right)^r = \left( \frac{\big(1+\lambda \log(1+t)\big)^{\frac{1}{\lambda }}-1}{2}+1\right)^{-r}\\
&= \sum_{m=0}^\infty {-r \choose m} 2^{-m} \left( \big(1+\lambda \log(1+t)\big)^{\frac{1}{\lambda }}-1 \right)^m\\
&= \sum_{m=0}^\infty (-1)^m 2^{-m} {r+m-1 \choose m} m! \sum_{k=m}^\infty S_{2,\lambda }(k,m) \frac{1}{k!} \big(\log(1+t)\big)^k\\
&=\sum_{k=0}^\infty \left( \sum_{m=0}^k (-1)^m m! {r+m-1 \choose m} 2^{-m} S_{2,\lambda }(k,m) \right) \sum_{n=k}^\infty S_1(n,k) \frac{t^n}{n!}\\
&=\sum_{n=0}^\infty \left\{ \sum_{k=0}^n \sum_{m=0}^k (-1)^m m! {r+m-1 \choose m} 2^{-m} S_{2,\lambda }(k,m) S_1(n,k) \right\} \frac{t^n}{n!}.
\end{split}\end{equation}
and

\begin{equation}\begin{split}\label{49}
&\left( \frac{2}{\big(1+\lambda \log(1+t)\big)^{\frac{1}{\lambda }}+1}\right)^r = \sum_{n=0}^\infty Ch_{n,\lambda }^{(r)} \frac{t^n}{n!}. 
\end{split}\end{equation}
Therefore, by \eqref{48} and \eqref{49}, we obtain the following theorem.

\begin{thm}
For $n \geq 0$, $r \in \mathbb{N}$, we have
\begin{equation*}\begin{split}
Ch_{n,\lambda }^{(r)} = \sum_{k=0}^n \sum_{m=0}^k (-1)^m m! {r+m-1 \choose m} 2^{-m} S_{2,\lambda }(k,m) S_1(n,k).
\end{split}\end{equation*}
\end{thm}

Let $\mathcal{E}_{m,\lambda }^{(r)}$ be the higher-order degenerate Euler numbers  defined by $\mathcal{E}_{n,\lambda }^{(r)} = \mathcal{E}_{n,\lambda }^{(r)}(0)$, $(n \geq 0)$.

Then, by \eqref{39}, we get

\begin{equation}\begin{split}\label{50}
&\sum_{n=0}^\infty \mathcal{E}_{n,\lambda }^{(r)} \frac{t^n}{n!} =
\left( \frac{2}{(1+\lambda t)^{\frac{1}{\lambda }}+1}\right)^r = \left( \frac{(1+\lambda t)^{\frac{1}{\lambda }}-1}{2}+1\right)^{-r}\\
&= \sum_{m=0}^\infty {-r \choose m} \big((1+\lambda t)^{\frac{1}{\lambda }}-1\big)^m 2^{-m}\\
&= \sum_{m=0}^\infty (-1)^m {r+m-1 \choose m} 2^{-m} m! \sum_{n=m}^\infty S_{2,\lambda }(n,m) \frac{t^n}{n!}\\
&= \sum_{n=0}^\infty \left( \sum_{m=0}^n (-1)^m m! {r+m-1 \choose m} 2^{-m} S_{2,\lambda }(n,m) \right) \frac{t^n}{n!}.
\end{split}\end{equation}

Thus, by \eqref{50}, we get

\begin{equation}\begin{split}\label{51}
\mathcal{E}_{n,\lambda }^{(r)} = \sum_{m=0}^n (-1)^m m! {r+m-1 \choose m} 2^{-m} S_{2,\lambda }(n,m).
\end{split}\end{equation}

From \eqref{37}, we have

\begin{equation}\begin{split}\label{52}
&\sum_{n=0}^\infty Ch_{n,\lambda }^{(r)}(x) \frac{t^n}{n!} = \left(\frac{2}{1+\big(1+\lambda \log(1+t)\big)^{\frac{1}{\lambda }}}\right)^r \big(1+\lambda \log(1+t)\big)^{\frac{x}{\lambda }} \\
&=\left(\frac{2}{1+\big(1+\lambda \log(1+t)\big)^{\frac{1}{\lambda }}}\right)^k \left(\frac{2}{1+\big(1+\lambda \log(1+t)\big)^{\frac{1}{\lambda }}}\right)^{r-k} \big(1+\lambda \log(1+t)\big)^{\frac{x}{\lambda }} \\
&= \left( \sum_{l=0}^\infty Ch_{l,\lambda }^{(k)} \frac{t^l}{l!} \right) 
\left( \sum_{m=0}^\infty Ch_{m,\lambda }^{(r-k)}(x) \frac{t^m}{m!} \right) \\
&= \sum_{n=0}^\infty \left( \sum_{l=0}^n {n \choose l} Ch_{l,\lambda }^{(k)}Ch_{n-l,\lambda }^{(r-k)}(x) \right) \frac{t^n}{n!}.
\end{split}\end{equation}

Therefore, by \eqref{52}, we obtain the following convolution result.

\begin{thm}

For $n \geq 0$, $r \in \mathbb{N}$, we have
\begin{equation*}\begin{split}
Ch_{n,\lambda }^{(r)}(x) =  \sum_{l=0}^n {n \choose l} Ch_{l,\lambda }^{(k)} Ch_{n-l,\lambda }^{(r-k)}(x).
\end{split}\end{equation*}
\end{thm}

Remark. By \eqref{37}, we easily get

\begin{equation}\begin{split}\label{53}
&\sum_{n=0}^\infty Ch_{n,\lambda }^{(r)}(x) \frac{t^n}{n!} = \left(\frac{2}{1+\big(1+\lambda \log(1+t)\big)^{\frac{1}{\lambda }}}\right)^r \big(1+\lambda \log(1+t)\big)^{\frac{x}{\lambda }} \\
&= \left( \sum_{l=0}^\infty Ch_{l,\lambda }^{(r)} \frac{t^l}{l!} \right) 
\left( \sum_{m=0}^\infty (x)_{m,\lambda } \frac{1}{m!} \big( \log(1+t)\big)^m \right)\\
&= \left( \sum_{l=0}^\infty Ch_{l,\lambda }^{(r)} \frac{t^l}{l!} \right) 
\left( \sum_{k=0}^\infty \left( \sum_{m=0}^k (x)_{m,k\lambda } S_1(k,m) \right) \frac{t^k}{k!} \right)\\
&= \sum_{n=0}^\infty \left( \sum_{k=0}^n \sum_{m=0}^k {n \choose k} (x)_{m,\lambda } S_1(k,m) Ch_{n-k,\lambda }^{(r)} \right) \frac{t^n}{n!}.
\end{split}\end{equation}

Comparing the coefficients on both sides of \eqref{53}, we have
\begin{equation*}\begin{split}
Ch_{n,\lambda }^{(r)}(x) =\sum_{k=0}^n \sum_{m=0}^k {n \choose k} (x)_{m,\lambda } S_1(k,m) Ch_{n-k,\lambda }^{(r)}.
\end{split}\end{equation*}


\begin{thebibliography}{0}

\bibitem{key-01}S. Araci, O. Ozer,
\textit{Extended $q$-Dedekind-type Daehee-Changhee sums associated with extended $q$-Euler polynomials.
} Adv. Difference Equ.
(2015) 2015:272 5 pp.

\bibitem{key-02} L. Carlitz,
\textit{Degenerate Stirling, Bernoulli and Eulerian numbers,
} Utilitas Math.
{\bf{15}} (1979), 51-88.


\bibitem{key-03} L. Carlitz,
\textit{A degenerate Staudt-Clausen theorem,
} Arch. Math. (Basel)
{\bf{7}} (1956), 28-33.

\bibitem{key-04} B. S. El-Desouky, A. Mustafa,
\textit{New Results on higher-order Daehee and Bernoulli numbers and polynomials,
} Adv. Difference Equ.
(2016) 2016:32, 21 pp.

\bibitem{key-05}D.-D. Su, Y. He,
\textit{Some new formulas for the products of the Frobenius-Euler polynomials,
} Adv. Difference Equ.
(2017) 2017:171.


\bibitem{key-06}B. M. Kim, J. Jeong, S.-H. Rim,
\textit{Some explicit identities on Changhee-Genocchi polynomials and numbers,
} Adv. Difference Equ.
(2016) 2016:202, 12 pp.

\bibitem{key-07}D. S. Kim, T. Kim, J.-J. Seo,
\textit{Note on Changhee polynomials and numbers,
} Adv. Studies Theor. Phys.
{\bf{7}} (2013), no. 20, 993-1003.


\bibitem{key-08}D. S. Kim, T. Kim,
\textit{Identities of symmetry for the generalized degenerate Euler polynomials,
} Modern Mathematical methods and high performance computing in science and technology,
35-43, Springer Proc. Math. Stat.,
{\bf{171}}, Springer, Singapore, 2016.

\bibitem{key-09}T. Kim, D. S. Kim,
\textit{A note on nonlinear Changhee differential equations,
}Russ. J. Math. Phys. 
{\bf{23}} (2016), no. 1, 88-92.

\bibitem{key-10} T. Kim,
\textit{A note on degenerate Stirling polynomials of the second kind,
} Proc. Jangjeon Math. Soc.
{\bf{20}} (2017), no. 3 (in press).

\bibitem{key-11}H.-I. Kwon, T. Kim, J.-J. Seo,
\textit{A note on degenerate Changhee numbers and polynomials,
} Proc. Jangjeon Math. Soc.
{\bf{18}} (2015), no. 3, 295-305.

\bibitem{key-12}J. Kwon, J.-W. Park,
\textit{On modified degenerate Changhee polynomials and numbers,
} J. Nonlinear Sci. Appl.
{\bf{9}} (2016), no. 12, 6294-6301.

\bibitem{key-13}J. G. Lee, L.-C. Jang, J.-J. Seo, S.-K. Choi, H.-I. Kwon,
\textit{On Appell-type Changhee polynomials and numbers,
} Adv. Difference Equ. 
(2016) 2016:160, 10 pp.

\bibitem{key-14} D. Lim,
\textit{Fourier series of higher-order Daehee and Changhee functions and their applications,
} J. Inequal. Appl. 
(2017), 2017:150.

\bibitem{key-15}J.-W. Park,
\textit{On the twisted $q$-Changhee polynomials of higher order,
} J. Comput. Anal. Appl.
{\bf{20}} (2016), no. 3, 424-431.

\bibitem{key-16}S.-H. Rim, J.-W. Park, S.-S. Pyo, J. Kwon,
\textit{The $n$-th twisted Changhee polynomials and numbers,
} Bull. Korean Math. Soc.
{\bf{52}} (2015), no. 3, pp. 741-749.

\bibitem{key-17}C. S. Ryoo, 
\textit{Symmetric identities for degenerate $(h,q)$-tangent polynomials associated with the $p$-adic integral on $\mathbb{Z}_p$,
} Inter. J. Math. Anal. 
{\bf{11}} (2017), no. 8, 353-362.

\bibitem{key-18} Y. Simsek,
\textit{Identities on the Changhee numbers and Apostol-type Daehee polynomials,
} Adv. Stud. Contemp. Math. (Kyungshang)
{\bf{27}} (2017), no. 2, 199-212.

\bibitem{key-19}Y. Simsek,
\textit{Analysis of the $p$-adic $q$-Volkenborn integrals: an approach to generalized Apostol-type special numbers and polynomials and their applications,
} Cogent Math.
{\bf{3}} (2016), Art. ID 1269393, 17 pp.

\bibitem{key-20} N. L. Wang, H. Li,
\textit{Some identities on the Higher-order Daehee and Changhee Numbers,
} Pure Appl. Math. J.
{\bf{4}} (2015), Issue 5-1, 33-37.


 
\end{thebibliography}
\end{document}